\documentclass[11pt]{article}
\usepackage{graphicx}
\usepackage{amssymb, amsmath, latexsym, amscd}
\newtheorem{question}{Question}

\def\G{\Gamma}
\def\Z{\mathbb Z}
\def\R{\mathbb R}
\def\Q{\mathbb Q}
\def\GL{${\rm{GL}}(n,{\Z})$}
\def\SL{${\rm{SL}}(n,{\Z})$}
\def\Out{${\rm{Out}}(F_n)$}
\def\Aut{${\rm{Aut}}(F_n)$}
\def\Map{${\rm{Mod}}(S_g)$}
\def\Gg{${\rm{Mod}}^{\pm }(S_g)$}

\def\F{B}

\title{Automorphism groups of free groups, surface groups and free abelian groups}
\author{Martin R. Bridson and Karen Vogtmann}
\date{23 July 2005}
\begin{document}
\maketitle

The group of $2 \times 2$ matrices with integer entries and determinant $\pm 1$ can be identified either with the group of outer automorphisms of a rank two free group or with the group of isotopy classes of homeomorphisms of a 2-dimensional torus.  Thus this group is the beginning of  three natural sequences of groups, namely the general linear groups ${\rm{GL}}(n,\Z)$, the groups \Out\  of outer automorphisms of free groups of rank $n\geq 2$, and the  mapping class groups \Gg\  of orientable surfaces of genus $g\geq 1$.
Much of the work on mapping class groups and automorphisms of free groups is motivated by the idea that these sequences of groups are strongly analogous, and should have many properties in common.  This program is occasionally derailed by uncooperative facts but  has in general proved to be a successful strategy, leading to fundamental discoveries about the structure of these groups.  In this article we will highlight a few of the most striking similarities and differences between these series of groups and present some open problems motivated by this philosophy.

Similarities among the groups  \Out, \GL\  and \Gg\  begin with the fact that these are the outer automorphism groups of the most primitive types  of torsion-free discrete groups, namely free groups, free abelian groups and the fundamental groups of closed orientable surfaces
$\pi_1S_g$. In the case of \Out\ and \GL\ this is obvious, in the case of
\Gg\  it is a classical theorem of Nielsen.
 In all cases there is a {\em{determinant}}
 homomorphism to ${\mathbb Z}/2$; the kernel of this map is the group of ``orientation-preserving"  or ``special"  automorphisms, 
and is denoted ${\rm{SOut}}(F_n),\, {\rm{SL}}(n,Z)$ or \Map\  respectively.   

\section{Geometric and topological models} 

A natural geometric  context for studying the global structure of
\GL\  is provided by the symmetric space $X$ of   positive-definite,
real symmetric matrices of determinant 1 (see \cite{Sou04} for a nice introduction to this subject). This is a 
non-positively curved manifold  diffeomorphic to   $\R^d$,
 where
$d={\frac 1 2 n(n+1)}-1$.  
\GL\  acts properly by isometries on $X$ with a quotient of finite volume.  

Each $A\in X$ defines an inner product on $\R^n$ and hence a  Riemannian metric
$\nu$ of constant curvature and 
volume 1 on the $n$-torus $T^n= \R^n/\Z^n$. One
can recover $A$ from the metric $\nu$ and an ordered basis for $\pi_1T^n$.
Thus $X$ is homeomorphic to the space of equivalence
classes of {\em marked} Euclidean tori $(T^n,\nu)$ of
volume 1, where a {\em marking} is a homotopy class of homeomorphisms
$\rho:T^n\to (T^n,\nu)$ and two marked tori are considered equivalent if there
is an isometry $i:(T_1^n,\nu_1)\to (T_2^n,\nu_2)$ such that $\rho_2^{-1}\circ i\circ\rho_1$
is homotopic to the identity.   The natural action of $\hbox{GL}(n,\Z)={\rm{Out}}(\Z^n)$
on $T^n=K(\Z^n,1)$  twists the markings on tori, and when one traces through
the identifications this is the standard action on $X$.

If one replaces $T^n$ by $S_g$ and
follows exactly this formalism with marked metrics of constant
curvature\footnote{if $g\ge 2$ then the curvature will be negative} and fixed volume,   then one
arrives at the definition of {\em Teichm\"uller space} and the natural action of 
\Gg$={\rm{Out}}(\pi_1S_g)$ on it. Teichm\"uller space is again homeomorphic to a Euclidean
space, this time $\R^{6g-6}$.

In the case of \Out\ there is no canonical choice of classifying space $K(F_n,1)$
but rather a finite collection of natural models, namely the finite graphs of 
genus $n$ with no vertices of valence less than 3. Nevertheless, one can proceed in
essentially the same way: one  considers metrics of fixed volume (sum of the
lengths of edges =1) on the various models for $K(F_n,1)$, each equipped
with a marking, and one makes the obvious identifications as the homeomorphism
type of a graph changes with a sequence of metrics that shrink an edge to
length zero. The space of marked metric structures obtained in this case is
Culler and Vogtmann's Outer space \cite{CulVog86}, which is stratified  by
manifold subspaces corresponding to the different homeomorphism types of
graphs that arise. This space is not   a manifold, but it  is  contractible and
its local homotopical structure is a natural generalization of that for a manifold (cf.~\cite{Vog90}).

One can also learn a great deal 
about the group \GL\ by
examining its actions on the Borel-Serre compactification of the symmetric space $X$ and on the spherical Tits building, which encodes the asymptotic geometry of $X$.
Teichm\"uller space and Outer space both admit useful bordifications that are closely analogous
to the Borel-Serre bordification \cite{Har88, Iva87, BesFei00}.  And in place of the spherical
Tits building for \GL\ one has the complex of curves
\cite{Har81} for  \Gg, which has played an important role in recent advances concerning the 
large scale geometry of \Gg. For the moment this complex has no well-established counterpart in the
context of \Out.

These closely parallel descriptions of geometries for the three families of groups
have led mathematicians to try to push the analogies further, both for the geometry and topology
of the ``symmetric spaces" and for   purely group-theoretic properties that are most
naturally proved using the geometry of the symmetric space.
For example, the symmetric space for \GL\  admits a natural equivariant deformation retraction onto an $n(n-1)/2$-dimensional cocompact subspace, the {\it well-rounded retract} \cite{Ash84}.  Similarly, both Outer space and the Teichm\"uller space of a punctured or bounded orientable surface retract equivariantly onto  cocompact simplicial spines \cite{CulVog86, Har88}.  In all these cases, the  retracts have dimension equal to the virtual cohomological dimension of the relevant group.    For closed surfaces, however, the question remains open:

\begin{question}\label{spine} 
Does the Teichm\"uller space for $S_g$
 admit an equivariant deformation retraction onto a cocompact 
spine whose dimension is equal to $4g-5$, the virtual cohomological
dimension of \Gg?
\end{question}
Further questions of a similar
nature are discussed in (\ref{BC}).  

The issues involved in using these symmetric space analogs to prove purely group theoretic
properties are illustrated in the proof of the Tits alternative, which holds for all three classes of groups.  A group $\Gamma$ is said to satisfy 
the Tits alternative if each of its subgroups
 either contains a non-abelian free 
 group or else is virtually solvable.  The strategy 
for proving this is similar in each of the three families that we
are considering:  inspired by Tits's
original proof for linear groups (such as \GL), one attempts to use a ping-pong argument on a suitable boundary at infinity of the symmetric space.  
This strategy ultimately succeeds but the details vary enormously between the
three contexts, and in the case of \Out\ they are particularly intricate 
 (\cite{BesFeiHan00, BesFeiHan97} versus \cite{BirLubMcC83}).
One finds that this is often the case:  analogies between the three classes of groups can   be carried through to theorems, and the architecture of the expected proof is often a good guide, but at a more detailed level the techniques required vary in essential ways from one class to the next and can be of completely different orders of difficulty.

Let us return to problems more directly phrased in terms of the geometry
of the symmetric spaces.
The symmetric space for \GL\  has a  left-invariant metric
of non-positive curvature, the geometry of which is relevant to many areas of mathematics beyond geometric group theory.  Teichm\"uller space has two natural metrics, the Teichm\"uller metric and
the Weyl-Petersen metric, and again  the study of each is a rich subject.  In contrast, the metric theory of Outer space has not been developed, and in fact there is no obvious candidate for a natural metric.  Thus, the following question has been left deliberately vague:

\begin{question}
Develop a metric theory of Outer space.
\end{question}

The elements of infinite order in \GL\ that are diagonalizable over $\mathbb C$
act as loxodromic isometries of $X$. When $n=2$, these elements
are the hyperbolic matrices; each fixes two points at infinity in 
$X=\mathbb H^2$,
one  a source and one a sink. The analogous type of element in \Gg\ is a
pseudo-Anosov, and in \Out\ it is an {\em iwip} (irreducible with irreducible powers).
In both cases, such elements have two fixed points at infinity (i.e.
in the natural boundary of the symmetric space analog), 
and the action of the
cyclic subgroup generated by the element exhibits  
the north-south dynamics
familiar from the action of hyperbolic matrices on the closure of
the Poincar\'e disc \cite{LevLus03}, \cite{Iva02}. 
In the case  of \Gg\  this cyclic subgroup leaves invariant a 
unique geodesic line in 
Teichm\"uller space, i.e. pseudo-Anosov's are axial like 
the semi-simple elements of infinite
order in \GL.
Work of Handel and Mosher \cite{HanMos05} shows that in the case
of {\em iwips} one cannot hope to have an axis in precisely the same metric sense,
but leaves open the possibility that   there may be a reasonable notion of  ``quasi-axis" for 
such automorphisms. 

\begin{question}  Find a useful description of a quasi-axis for an iwip
acting on Outer Space,  with limit set the fixed points of the iwip at
infinity.
\end{question}

\section{Actions of \Aut\ and \Out\ on other spaces}

Some of the questions that we shall present are more naturally stated
in  terms of \Aut\ than \Out, while some are natural for both. To
avoid redundancy, we shall state only one form of each question. 

\subsection{Baum-Connes and Novikov conjectures}\label{BC}

Two famous conjectures relating topology, geometry and functional analysis are the Novikov and Baum-Connes conjectures.  The Novikov conjecture  for closed oriented manifolds  with fundamental group $\Gamma$ says that certain {\it higher signatures\ } coming from $H^*(\Gamma;\Q)$ are homotopy invariants.  It is implied by the Baum-Connes conjecture, which says that a certain {\it assembly map\ } between two  $K$-theoretic objects associated 
to $\Gamma$  is an isomorphism. Kasparov \cite{Kasp88} 
proved the Novikov conjecture for \GL, and  Guenther, Higson and
Weinberger proved it for all linear groups \cite{GHW05}.
The Baum-Connes conjecture for \GL\ is open when $n\ge 4$ 
(cf.~\cite{Laff98}).

 Recently Storm \cite{Sto05} pointed out that  the Novikov 
conjecture for mapping class groups
follows from results that have been announced by Hamenst\"adt 
\cite{Ham04} and Kato \cite{Kat00}, leaving open the following: 

\begin{question} Do mapping class groups or \Out\ satisfy the Baum-Connes conjecture? 
Does \Out\ satisfy the Novikov conjecture?
\end{question}
An approach to proving these conjectures is given by work of  
Rosenthal  \cite{Ros05}, 
generalizing results
of Carlsson and Pedersen \cite{CarPed95}.  A contractible space on which a group $\G$ 
acts properly and for which the fixed point sets of finite subgroups are contractible is called an $\underbar{E}\Gamma$. 
Rosenthal's theorem says that the 
Baum-Connes map for $\G $
 is split injective if there is a cocompact $\underbar{E}\Gamma=E$ that admits  a compactification $X$, such that 
\begin{enumerate}
\item the $\G$-action extends to $X$;
\item $X$ is metrizable;
\item $X^G$ is contractible for every finite subgroup $G$ of $\Gamma$
\item $E^G$ is dense in $X^G$ for every finite subgroup $G$ of $\Gamma$
\item compact subsets of E become small near $Y = X\smallsetminus
E$ under the $\Gamma$-action: for every
compact $ K \subset E$ and every neighborhood $U\subset X$ of $y\in Y$, there
exists a neighborhood $V\subset X$ of $y$ such 
that  $\gamma K \cap V \neq\emptyset$ implies
$\gamma K \subset  U$.
\end{enumerate}
The existence of such a space $E$ also implies the Novikov conjecture for $\Gamma$.

For   \Out\ the spine of Outer space mentioned in the previous section is a reasonable candidate for the required  $\underbar {E}\G$, and there is a similarly defined candidate for \Aut.  For mapping class groups of punctured  surfaces  the complex of arc systems which fill up the surface is a good candidate (note that this can be identified with a subcomplex of Outer space, as in \cite{HatVog96}, section 5).

\begin{question} Does there exist
a compactification of the spine of Outer space satisfying Rosenthal's conditions? Same question for the complex of arc systems filling a punctured surface.
\end{question}

 In all of the cases mentioned above, the candidate space $E$  has dimension equal to the virtual cohomological dimension of the group.  G.~Mislin \cite{Mis04} has 
constructed a cocompact $\underbar{E}G$ for the mapping class group of a closed surface, but it has much higher dimension, equal to the dimension of the Teichm\"uller space.
This leads us to a slight variation on Question \ref{spine}.

\begin{question} Can one
 construct a cocompact $\underbar{E}G$ with dimension equal  to the 
virtual cohomological dimension of the mapping class group of a closed surface?
\end{question}

\subsection{Properties (T) and FA}

A group has Kazdhan's Property (T) if any action of the group 
by isometries on a Hilbert
space has fixed vectors. Kazdhan proved that \GL\ has property (T)
for $n\ge 3$. 

\begin{question} For $n>3$, does \Aut\  have property (T)?
\end{question}

The corresponding question for mapping class groups is also
open. If  \Aut\   were to
have Property (T), then an argument of Lubotzky and Pak \cite{LubPak01} would
provide a conceptual explanation
of the apparently-unreasonable effectiveness of certain 
algorithms in computer science, specifically 
 the Product Replacement Algorithm of Leedham-Green {\em{et al}}.

If a group has Property (T) then it has Serre's
property FA:  every action of the group on an $\mathbb R$-tree has a fixed point. 
When $n\ge 3$, \GL\  has property FA, as
do \Aut\ and \Out, and mapping class groups in genus $\ge 3$ (see \cite{CulVog96}).
In contrast,  McCool \cite{McC89} has shown
that ${\rm{Aut}}(F_3)$ has a subgroup of finite-index with positive first
betti number, i.e. a subgroup which
 maps onto $\mathbb Z$. In particular this subgroup acts by translations
on the line and therefore does not have property FA or (T). Since
property (T) passes to finite-index subgroups, it follows 
that ${\rm{Aut}}(F_3)$ does not have property (T).

\begin{question}\label{betti} 
 For $n>3$, 
does \Aut\ have a subgroup of finite index with positive first betti number?
\end{question}

Another finite-index subgroup of $Aut(F_3)$ mapping onto $\mathbb Z$ was constructed by Alex Lubotzky, and was explained to us by
Andrew Casson. Regard $F_3$ as the fundamental group of a graph $R$ with
one vertex. The single-edge loops 
 provide a basis
 $\{a,b,c\}$ for $F_3$. Consider the 2-sheeted
covering $\hat R\to R$ with fundamental group  
$\langle a,b,c^2,cac^{-1},cbc^{-1}
\rangle$ and let $G\subset{\rm{Aut}}(F_3)$ be the stabilizer of
this subgroup. $G$ acts on $H_1(\hat R,\mathbb Q)$
leaving invariant
the eigenspaces of the involution that generates the Galois group of the
covering. The eigenspace corresponding to the eigenvalue $-1$ is two
dimensional with basis $\{a-cac^{-1},\, b-cbc^{-1}\}$. The action of $G$
with respect to this basis gives  an epimorphism $G\to{\rm{GL}}(2,\Z)$.
Since ${\rm{GL}}(2,\Z)$ has a free subgroup of
finite-index, we  obtain a subgroup of finite index in ${\rm{Aut}}
(F_3)$ that
maps onto a non-abelian free group.

One can imitate  the essential features of this construction with various
other finite-index subgroups of $F_n$, thus producing subgroups of
finite index in \Aut\ that map onto ${\rm{GL}}(m,\Z)$. In each case one
finds  that $m\ge n-1$.

\begin{question} If there is a homomorphism from a 
subgroup of finite index in \Aut\ onto a subgroup of finite index
in $GL(m,\Z)$, then must $m\ge n-1$? 
\end{question}

Indeed one might ask:

\begin{question} If $m<n-1$ and $H\subset$\Aut\ is a subgroup of
finite index,  then does every  homomorphism 
 $H\to{\rm{GL}}(m,\Z)$ have finite image?
\end{question}

Similar questions are interesting for the other groups in our
families (cf.~section 3). For example, If $m<n-1$ and $H\subset$\Aut\ is a subgroup of
finite index,  then does every  homomorphism 
 $H\to{\rm{Aut}}(F_m)$ have finite image?

A positive answer to the
following question would answer Question \ref{betti}; a negative answer
would show that \Aut\ does not have property (T).

\begin{question} For $n\ge 4$,
do subgroups of finite index in \Aut\ have Property FA? 
\end{question}

 A promising approach to this last question breaks down because we do not know
the answer to the following question.

\begin{question} Fix a basis for $F_n$ and let $A_{n-1}\subset
Aut(F_n)$ be the copy of $Aut(F_{n-1})$ corresponding to the first
$n-1$ basis elements. Let $\phi: Aut(F_n)\to G$ be a homomorphism 
of groups. If $\phi(A_{n-1})$ is finite, must the image of $\phi$ be finite?
\end{question}

 Note that  the obvious analog of
this question for  \GL\ has a positive answer and plays a  role in the foundations of 
algebraic $K$-theory.

A different approach to establishing Property (T) was developed by Zuk \cite{Zuk96}.   He established a combinatorial criterion on the links of vertices in a simply connected $G$-complex which, if satisfied, implies
 that $G$ has property (T): one must show that the smallest positive eigenvalue of the discrete Laplacian on links is sufficiently large.  In addition to the $Aut(F_n)$ analog of the spine of Outer space, there are several other simply-connected  complexes
on which \Aut\ acts, and these  might be used to test Zuk's criterion.  
Hand-worked experiments enable one to get arbitrarily close to the critical value
in Zuk's criterion as $n\to\infty$, but it is not clear how to interpret this evidence.

\subsection{Actions on CAT$(0)$ spaces}

An $\mathbb R$-tree may be defined as a complete 
CAT$(0)$ space of dimension\footnote{topological covering dimension} 1. 
Thus one might generalize property FA by asking, for each
$d\in\mathbb N$, which groups must fix a point
whenever they act by isometries on a complete CAT$(0)$
space of dimension $\le d$.

\begin{question}
What is the least integer $\delta$ such that \Out\ acts without
a global fixed point on a complete CAT$(0)$ space of dimension $\delta$?
And what is the least dimension for the mapping class group \Gg?
\end{question}

The action of \Out\ on the first homology of $F_n$ defines a map
from ${\rm{Out}}(F_n)$ to ${\rm{GL}}(n,\mathbb Z)$ and hence an
action of \Out\ on the symmetric space for ${\rm{GL}}(n,\mathbb R)$,
which is a complete CAT$(0)$ space of dimension $\frac1 2{n(n+1)}-1$. This
action does not have a global fixed point and hence we obtain a
quadratic upper bound  $n(n+1)/2$ on $\delta$.  On the other hand, since \Out\ 
has property FA, $\delta\ge 2$. In fact, motivated by work of 
Farb on ${\rm{GL}}(n,\mathbb Z)$, Bridson \cite{Brid05} has
shown that using a Helly-type theorem and the structure
of finite subgroups in \Out, one can obtain a lower bound on $\delta$ that
grows as a linear function of $n$.
 Note that a lower bound of $3n-3$ on $\delta$
would imply that Outer Space did not support a complete
\Out-equivariant metric of non-positive
curvature. 

If $X$ is a CAT$(0)$ polyhedral complex with only finitely many isometry types of cells 
(e.g. a finite dimensional cube complex), then each isometry of $X$
is either elliptic (fixes a point) or hyperbolic (has an axis of translation) \cite{Bri99}.
If $n\ge 4$ then a variation on an argument of Gersten \cite{Ger94} shows that
in any action of \Out\ on $X$, no
Nielsen generator can act as a hyperbolic isometry. 

\begin{question} If $n\ge 4$, then can \Out\ act without a global
fixed point on a finite-dimensional  {\rm{CAT$(0)$}}  cube complex?
\end{question}

\subsection{Linearity}

 Formanek and Procesi \cite{ForPro92}  proved that \Aut\ is not linear for $n\geq 3$
by showing that ${\rm{Aut}}(F_3)$ 
contains a ``poison subgroup", i.e. a subgroup which has no faithful linear representation.  

 Since ${\rm{Aut}}(F_n)$ embeds in ${\rm{Out}}(F_{n+1})$, this settles 
the question of linearity for 
\Out\ as well, except when $n=3$.

\begin{question} Does $\rm{Out}(F_3)$ have a faithful representation
into $\rm{GL}(m,\mathbb C)$ for some $m\in\mathbb N$?
\end{question}

Note that braid groups are linear \cite{Big01} but it is unknown if mapping
class groups of closed surfaces are.  Brendle and Hamidi-Tehrani  \cite{BreHam01}
 showed that the approach of Formanek and Procesi cannot be adapted
directly to the mapping class groups. More precisely,
they prove that  the type of ``poison subgroup" described above does not
arise in   mapping class groups.

The fact that the above question  remains open is an indication that
${\rm{Out}}(F_3)$ can behave differently from \Out\ for $n$ large;
the existence of finite index subgroups mapping onto $\Z$ was
another instance of this, and we shall see another in our discussion of
automatic structures and isoperimetric inequalities.

\section{Maps to and from \Out}

A particularly intriguing aspect of the analogy between \GL\  and the
two other classes of groups is the extent to which the celebrated
rigidity phenomena for lattices in higher rank semisimple groups
transfer to mapping class groups and \Out. Many of the questions
in this section concern aspects of this rigidity; questions 9 to 11
should also be viewed in this light.
 
Bridson and Vogtmann \cite{BriVog03:2} showed that any homomorphism
from \Aut\ to a group $G$ has finite image if $G$ does not contain the symmetric group $\Sigma_{n+1}$; in particular, any homomorphism ${\rm{Aut}}(F_n)\to
{\rm{Aut}}(F_{n-1})$ has image of order at most 2.

\begin{question} If $n\ge 4$ and $g\ge 1$, does every homomorphism from
 \Aut\ to \Gg\  have finite image?
\end{question}

By \cite{BriVog03:2}, one cannot obtain homomorphisms with infinite
image unless \Gg\   contains the symmetric group $\Sigma_{n+1}$.  For large enough genus, you can realize any symmetric group; but the order of a finite group of symmetries is at most 84g-6, so here
one needs $84g-6\geq (n+1)!$.

There are no {\em injective} maps from ${\rm{Aut}}(F_n)$ to
 mapping class groups.  This follows from the result of 
Brendle and Hamidi-Tehrani  that we quoted  earlier.
For certain $g$ one can construct homomorphisms 
${\rm{Aut}}(F_3)\to$\Gg\  with infinite image, but we do not know the
minimal such $g$.

 \begin{question} Let $\Gamma$ be an irreducible lattice in a semisimple Lie group of $\R$-rank at least 2. Does every homomorphism from $\Gamma$
to \Out\ have finite image?
\end{question}

This is known for non-uniform lattices (see \cite{BriFar01}; it follows easily
from the Kazdhan-Margulis finiteness theorem and the fact that solvable
subgroups of \Out\ are virtually abelian \cite{BesFeiHan04}).  
Farb and Masur  provided a positive answer to the analogous 
question for maps to mapping
class groups \cite{FarMas98}. The proof of their theorem was based on
results of Kaimanovich and Masur \cite{KaiMas96}
concerning random walks on Teichm\"uller
space.  (See  \cite{Iva02} and, for an alternative
approach, \cite{BesFuj02}.)

\begin{question} Is there a theory of random walks on Outer space similar to that of Kaimanovich and Masur for Teichm\"uller space?
\end{question}

Perhaps the most promising approach to Question 17 is
via bounded cohomology, following the template of Bestvina
and Fujiwara's work on subgroups of the mapping class
group \cite{BesFuj02}.

\begin{question} If a subgroup $G\subset$\Out\ is not
virtually abelian, then is $H^2_b(G;\mathbb R)$ infinite
dimensional?
\end{question}

If $m\ge n$ then
there are obvious embeddings ${\rm{GL}}(n,\Z)\to {\rm{GL}}(m,\Z)$
and ${\rm{Aut}}(F_n)\to{\rm{Aut}}(F_m)$, but there are no obvious embeddings
\Out\ $\to{\rm{Out}}(F_m)$.  Bogopolski and Puga \cite{BogPug04}
have shown that, for $m = 1 + (n-1)kn$, where $k$
 is an arbitrary natural number coprime to $n-1$, there is in fact
an embedding, by restricting automorphisms to a suitable characteristic subgroup of $F_m$.

\begin{question} For which values of $m$ does
 \Out\ embed in ${\rm{Out}}(F_m)$?  What is the
minimal such $m$, and is it true for all sufficiently large $m$?
\end{question}

Hatcher and Vogtmann \cite{HatVog04} showed that when
$n$ sufficiently large with respect to $i$,  the
 homology group $H_i({\rm{Out}}(F_n),\Z)$
is independent of $n$.
 
 \begin{question} Is there a map
 \Out\ $ \to{\rm{Out}}(F_m)$ that induces an isomorphism on homology
in the stable range?  
\end{question}

A number of the questions in this section and (2.2) ask whether certain quotients of
\Out\ or \Aut\ are necessarily finite. The following quotients arise naturally
in this setting: define $Q(n,m)$ to be the quotient of \Aut\ by the normal
closure of $\lambda^m$, where $\lambda$ is the Nielsen move defined on
a basis $\{a_1,\dots,a_n\}$ by $a_1\mapsto a_2a_1$. (All such Nielsen moves
are conjugate in \Aut, so the choice of basis does not alter the quotient.)

The image of a Nielsen move in \GL\ is an elementary matrix and the quotient
of \GL\ by the normal subgroup generated by the $m$-th powers of the
elementary matrices is the finite group ${\rm{GL}}(n,\mathbb Z/m)$. But
Bridson and Vogtmann \cite{BriVog03:2}
showed  that if $m$ is sufficiently large then $Q(n,m)$ is infinite
because it has a quotient that contains a copy of the free Burnside group 
$B(n-1,m)$.  Some further information can be gained by replacing $B(n-1,m)$
with the quotients of $F_n$ considered in subsection 39.3 of 
A.Yu.~Ol'shanskii's book \cite{Ols91}. But we   know very little about 
the groups $Q(n,m)$. For example:
 
\begin{question} For which values of $n$ and $m$ is $Q(n,m)$ infinite?
Is $Q(3,5)$ infinite?
\end{question}

\begin{question} Can $Q(n,m)$ have infinitely
many finite quotients? Is it residually finite? 
\end{question}

\section{Individual elements and mapping tori}\label{growth}

Individual elements $\alpha\in{\rm{GL}}(n,\mathbb Z)$ can be realized
as diffeomorphisms $\hat\alpha$ of the $n$-torus, while individual
elements $\psi\in$ \Gg\  can be realized as diffeomorphisms $\hat\psi$ of the surface
$S_g$. Thus one can study $\alpha$ via the geometry of the torus
bundle over $\mathbb S^1$ with holonomy $\hat\alpha$ and one can
study $\psi$ via the geometry of the 3-manifold that fibres over 
$\mathbb S^1$ with holonomy $\hat\psi$. (In each case the manifold
depends only on the conjugacy class of the element.)

The situation for \Aut\ and \Out\ 
is more complicated: the natural choices of classifying space $Y=K(F_n,1)$ 
are finite graphs of genus $n$, and no element of infinite order $\phi\in$\Out\ 
is induced by the action on $\pi_1(Y)$ of a 
homeomorphism of $Y$. Thus the best that
one can hope for in this situation is to identify a graph $Y_\phi$ that admits a 
homotopy equivalence inducing $\phi$ and that has additional structure
well-adapted to $\phi$. One would then form the mapping
torus of this homotopy equivalence to get a good classifying space for the
algebraic mapping torus $F_n\rtimes_\phi\mathbb Z$.

The {\em train track technology} of
 Bestvina, Feighn and Handel \cite{BesHan92, BesFeiHan00, BesFeiHan97}  is a major piece of work that derives
suitable graphs $Y_\phi$ with additional structure
encoding key properties of $\phi$. This results in a decomposition
theory for elements of \Out\ that is closely analogous to (but 
more complicated than) the Nielsen-Thurston theory for surface
automorphisms. Many of the results mentioned  in this section are premised on a 
detailed knowledge of this technology and one expects that a resolution
of the questions will be too.

There are several natural ways to
define the {\em growth} of an automorphism $\phi$ of a group $G$ with
finite generating set $A$; in the case of free, free-abelian, and surface groups 
these are all asymptotically
equivalent. The most easily defined growth function is
 $\gamma_\phi(k)$ where
$\gamma_\phi(k):=\max\{d(1,\phi^k(a)\mid a\in A\}$. 
 If $G=\Z^n$
then $\gamma_\phi(k)\simeq k^d$ for some integer $d\le n-1$, or else
$\gamma_\phi(k)$ grows exponentially. If $G$ is a surface group, the
Nielsen-Thurston theory shows that
only bounded, linear and exponential growth can occur. If $G=F_n$ 
and $\phi\in{\rm{Aut}}(F_n)$ then, as in the abelian case,  
$\gamma_\phi(k)\simeq k^d$ for some integer $d\le n-1$ or else
$\gamma_\phi(k)$ grows exponentially. 

\begin{question} Can one detect the growth of a surface or free-group
homomorphism by its action on the homology of a characteristic subgroup
of finite index?
\end{question}

Notice that one has to pass to a subgroup of finite index in order to
have any hope because  
 automorphisms of   exponential growth can act trivially on homology.
A.~Piggott \cite{Pig04} has answered the above question for free-group automorphisms
of polynomial growth, and linear-growth automorphisms
of surfaces are easily dealt with, but the exponential case remains open
in both settings. 

Finer questions concerning growth are addressed in the on-going work of
Handel and Mosher \cite{HanMos05}. They explore, for example, the implications of the
following contrast in behaviour between surface automorphisms and 
free-group automorphisms: in the surface case the exponential growth
rate of a pseudo-Anosov automorphism is the same as that of its
inverse, but this is not the case for iwip free-group automorphisms.

For mapping tori of automorphisms of free abelian groups $G=\Z^n\rtimes_\phi\Z$,
the following conditions are equivalent (see \cite{BriGer96}): $G$ is automatic; $G$  is
a CAT$(0)$ group\footnote{this means that $G$ acts properly 
and cocompactly by
isometries on a CAT$(0)$ space}; $G$ satisfies a quadratic isoperimetric inequality.
In the case of mapping tori of surface automorphisms, all mapping tori
satisfy the first and last of these conditions and one understands exactly
which $S_g\rtimes\mathbb Z$ are CAT$(0)$
groups. 

Brady, Bridson and Reeves  \cite{BraBriRee05} 
show that there exist mapping tori of free-group
automorphisms $F\rtimes\mathbb Z$ that are not automatic, and Gersten
showed that some are not CAT$(0)$ groups \cite{Ger94}. On the other hand, many
such groups do have these properties, and we have already noted that
they all satisfy a quadratic isoperimetric inequality.

\begin{question} Classify those $\phi\in{\rm{Aut}}(F_n)$ for which
$F_n\rtimes_\phi\Z$ is automatic and those for which it is {\rm{CAT}}$(0)$.
\end{question}

Of central importance in trying to understand mapping tori is:

\begin{question} Is there an alogrithm to decide isomorphism among groups
of the form $F\rtimes\mathbb Z$.
\end{question}

In the purest form of this question one is given the groups as finite
presentations, so one has to address issues of how to find the decomposition
$F\rtimes\mathbb Z$ and one has to combat the fact that this decomposition
may not be unique. But the heart of any solution should be an answer
to:

\begin{question} Is the conjugacy problem solvable in \Out?
\end{question}

Martin Lustig posted a detailed outline of a solution to this problem 
on his web page some years ago \cite{Lus99}, but neither this proof nor any other has been accepted for publication.  This problem is of central importance to the field and a clear, compelling solution would be of great interest.  
The conjugacy problem for mapping class groups
was shown to be solvable by Hemion \cite{Hem79}, and an effective algorithm for determining conjugacy, at least for pseudo-Anosov mapping classes, was given by Mosher \cite{Mos84}.  The isomorphism problem for groups of the
form $S_g\rtimes\mathbb Z$ can be viewed as  a particular case of the
solution to the isomorphism problem for fundamental groups of geometrizable
3-manifolds \cite{Sel95}.
 The solvability of the conjugacy problem for {\rm{GL}}$(n,\mathbb Z)$
and of the isomorphism problem among groups of the form $\Z^n\rtimes\Z$
is classical.

\section{Cohomology} 

In each of the series of groups ${\{\G_n\}}$  we are considering, 
the $i$th homology of $\G_n$ has been shown to be independent of $n$ for $n$ sufficiently large. 
For \GL\ this is due to Charney \cite{Cha79}, for mapping class groups to Harer \cite{Har90}, and for \Aut\ and \Out\  to
Hatcher and Vogtmann \cite{HatVog98:3, HatVog04}.   With trivial rational coefficients, 
the stable cohomology of \GL\ was computed in the 1970's by Borel \cite{Bor74},  and the stable 
rational cohomology of the mapping class group computed by Madsen and Weiss in 2002 \cite{MadWei04}.  The question for \Out\ remains open:

\begin{question}
What is the stable rational cohomology of \Out?
\end{question}
No non-trivial stable rational cohomology classes have been found, and the standard conjecture is that
the stable rational cohomology is in fact trivial.  

The exact stable range for trivial rational coefficients is known for \GL\ and for mapping class groups of punctured surfaces.
For \Out\ the best known result is that the $i$th homology is independent of $n$ for  $n>5i/4$, but
the exact range is unknown:

\begin{question}
Where precisely does the rational homology of \Out\ stabilize? And for \Aut?
\end{question}

Since the stable rational cohomology of mapping class groups and \GL\ is non-trivial, there is a natural impulse to use these maps to try to detect cohomology in \Out.  However, a result of Igusa's \cite{Igu02} shows that this does not work for \GL:  the map from \Out\ to \GL\ is trivial on homology in the stable range.  Wahl \cite{Wah04} has studied the 
mapping class group case
and shown that the map is an infinite loop space map. But the following question remains open.

\begin{question}  Do the natural maps from the  mapping class groups
of non-closed surfaces to  \Out\  induce the zero map on stable rational homology?
\end{question}

In fact, there are only two known non-trivial classes in the rational homology of \Out\  \cite{HatVog98:2, ConVog04}, both below the stable range.   
However, Morita \cite{Mor99} has defined an infinite series of cycles, by  using work of Kontsevich which identifies the cohomology of \Out\ with the homology of a 
certain infinite-dimensional Lie algebra.  The first of these cycles coincides with the only previously known cohomology class, in $H_4(
{\rm{Out}}(F_4);\Q)$, and Conant and Vogtmann \cite{ConVog04} showed that the second also gives a non-trivial class, in $H_8({\rm{Out}}(F_6);\Q)$.  Both Morita and Conant-Vogtmann also defined more general cycles, parametrized by odd-valent graphs.  

\begin{question} Are Morita's original cycles non-trivial in homology? Are the generalizations due to Morita and to Conant and Vogtmann non-trivial in homology?
\end{question} 

We note that Morita has identified several conjectural relationships between his cycles and various other interesting objects, including the image of the Johnson homomorphism, the group of homology cobordism classes of homology cylinders,  and the motivic Lie algebra associated to the algebraic mapping class group (see Morita's article in this volume).

It is interesting to note that $H_8({\rm{GL}}(6,\Z);\Q)\cong
 \Q$ \cite{ElbGanSou02}; this leads naturally to the question

\begin{question}
Is the image of the second Morita class in $H_8({\rm{GL}}(6,\Z;\Q))$ non-trivial? 
\end{question}

Finally, it is worth noting that some of Conant and Vogtmann's generalizations of the Morita cycles lie in the stable range.

\section{Generators and Relations}

The groups we are considering are all finitely generated. In each case, the most natural set of generators consists of a single orientation-reversing generator of order two,  together with a collection of simple infinite-order special automorphisms.  For \Out, these special automorphisms are the Nielsen automorphisms, which multiply one generator of $F_n$ by another and leave the rest of the generators fixed;  for \GL\  these are the elementary matrices; and for mapping class groups they are Dehn twists around a small set of 
non-separating simple closed curves.  

These generating sets have a number of
important features in common. First, implicit in the description of each is
a choice of generating set for the group $\F$ on which $\G$ is acting. In the case 
of \Gg\  this ``basis" can be taken to consist of $2g+1$ simple closed curves
representing  the standard generators $a_1,b_1,a_2,b_2, \ldots,a_g,b_g,$
of $\pi_1(S_g)$ together with $z=a_2^{-1}b_{3}a_{3}b_{3}^{-1}$.
In the case of \Out\  and 
\GL, the generating set is a basis for $F_n$ and $\Z^n$ respectively.

Note that in the cases $\G=$\Out\  or 
\GL, the universal property of the underlying free objects $\F=F_n$ or $\Z^n$
ensures that $\G$ acts transitively on the set of preferred generating sets (bases).
In the case $\F=\pi_1S_g$, the corresponding result is that any two collections of
simple closed curves with the same pattern of intersection numbers and 
 complementary regions are related by a homeomorphism of the surface, 
hence (at the level of $\pi_1$) by the action of $\G$.

If we identify $\Z^n$ with  the abelianization of $F_n$ and choose bases
accordingly, then the action of \Out\ on the abelianization induces a 
homomorphism $\hbox{\Out}\to\hbox{\GL}$ that sends each Nielsen move to
the corresponding elementary matrix (and hence is surjective). Correspondingly,
the action \Gg\  on the abelianization of $\pi_1S_g$ yields a homomorphism
onto the symplectic group $Sp(2g,\Z)$ sending the generators  of \Gg\  given by Dehn twists around the $a_i$ and $b_i$ to
transvections.
Another common feature of these generating sets is that they all have linear
growth (see  section \ref{growth}).

Smaller (but less transparent) generating sets exist in each case. 
Indeed 
B.H.~Neumann \cite{Neu32}
proved that \Aut\ (hence its quotients \Out\  and \GL) is generated by
just 2 elements  when $n\ge 4$. 
Wajnryb \cite{Waj96} proved that this is also
true of mapping class groups.

In each case one can also find generating 
sets consisting of finite order elements,
and in fact by involutions.  Zucca showed that \Aut\ can be generated by 3 involutions two of which commute \cite{Zuc97},
and Kassabov, building on work of Farb and Brendle, showed that mapping class groups of large enough  genus can be generated by 4 involutions \cite{Kas03}.

Our groups are also all finitely presented.  For \GL, or more precisely for \SL,  there are the classical Steinberg relations, which involve commutators of the elementary matrices.  For the special automorphisms ${\rm{SAut}}(F_n)$, Gersten gave a presentation in terms of corresponding commutator relations of the Nielsen generators \cite{Ger84}.  Finite presentations of the mapping class groups are more complicated. The first was given by Hatcher and Thurston,  and worked out explicitly by Wajnryb \cite{Waj83}.  

\begin{question} Is  there a set of simple Steinberg-type relations for the mapping class group?    
\end{question}

There is also a presentation of \Aut\ coming from the action of \Aut\  on
the subcomplex of Auter space spanned by graphs of degree at most 2.  This is simply-connected by \cite{HatVog98:3},  so Brown's method \cite{Bro84} can be used to write down a presentation. The vertex groups are stabilizers of marked graphs, and the edge groups are the stabilizers of pairs consisting of a marked graph and a forest in the graph.  The quotient of the subcomplex modulo \Aut\ can be computed explicitly, and one finds that  \Aut\ is generated by the (finite) stabilizers of seven specific marked graphs.  In addition, all of the relations except two come from the natural inclusions of edge stabilizers into vertex stabilizers, i.e. either including the stabilizer of a pair (graph, forest) into the stabilizer of the graph, or into the stabilizer of the quotient of the graph modulo the forest.  
Thus  the whole group is almost (but not quite) a pushout of these finite subgroups. In the terminology of Haefliger (see \cite{BriHae99},  II.12), 
the complex of groups is not simple.

\begin{question} Can \Out\ and \Gg\ be obtained as a pushout  of a finite subsystem
of their finite subgroups, i.e. is either the fundamental group of a developable simple
complex of finite groups on a 1-connected base?
\end{question}

\subsection{IA automorphisms}

We conclude with a well-known problem about the kernel ${\rm{IA}}(n)$
 of the map from ${\rm{Out}}(F_n)$  to ${\rm{GL}}(n,Z)$.  The notation ``IA" stands for {\it identity on the abelianization}; these are (outer) automorphisms of $F_n$ which are the identity on the abelianization $Z^n$ of $F_n$.
 Magnus showed that this 
 kernel is finitely generated, and for $n=3$ Krstic and 
 McCool showed that it is not finitely presentable 
 \cite{KrsMcC97}.  It is also known that in some dimension 
 the homology is not finitely generated \cite{SmiVog87}.
  But that is the extent of our knowledge of basic finiteness properties.

\begin{question}
Establish finiteness properties of  the kernel ${\rm{IA}}(n)$ 
of the map from ${\rm{Out}}(F_n)$ to ${\rm{GL}}(n,\Z)$.  In particular,
 determine whether ${\rm{IA}}(n)$ is finitely presentable for $n>3$.
\end{question}  

The subgroup ${\rm{IA}}(n)$ is analogous to the 
Torelli subgroup of the mapping class group of a surface, 
which also remains quite mysterious in spite of having
 been extensively studied.

\section{Automaticity and Isoperimetric Inequalities}

In the foundational text on automatic groups \cite{EpEtAl92},
Epstein gives a detailed account of Thurston's proof 
 that if $n\ge 3$ then \GL\ is not automatic. 
The argument uses
 the geometry of the symmetric space to obtain an
exponential lower bound on the $(n-1)$-dimensional isoperimetric
function of \GL; in particular the Dehn function of 
${\rm{GL}}(3,\mathbb Z)$ is shown to be exponential.

Bridson and Vogtmann \cite{BriVog95}, building on
this last result, proved that the Dehn functions of 
${\rm{Aut}}(F_3)$ and ${\rm{Out}}(F_3)$ are exponential.
They also proved that for all $n\ge 3$, 
neither \Aut\ nor \Out\ is biautomatic. 
In contrast, Mosher proved that mapping class groups
are automatic \cite{Mosh95} and Hamenst\"adt \cite{Ham04}
proved that they are biautomatic; in particular these
groups have quadratic Dehn functions and satisfy a polynomial
isoperimetric inequality in every dimension. Hatcher
and Vogtmann \cite{HatVog96} obtain an exponential upper
bound on the isoperimetric function of \Aut\ and \Out\ in
every dimension. 

An argument  sketched by Thurston
and expanded upon by Gromov \cite{Gromov93},
\cite{Gromov00}  (cf.~\cite{Drutu04})  indicates that
the Dehn function of \GL\ is quadratic when $n\ge 4$.
More generally,  the isoperimetric functions
of \GL\ should parallel those of Euclidean space 
in dimensions $m\le n/2$.

\begin{question}
What are the Dehn functions of \Aut\ and \Out\ for $n>3$?
\end{question}

\begin{question} 
 What are the higher-dimensional isoperimetric functions of
\GL, \Aut and \Out?
 \end{question}

\begin{question}  Is \Aut\ automatic for $n>3$?
\end{question}

\bibliographystyle{siam}
\bibliography{bibliography}

MRB: Mathematics.
Huxley Building,
Imperial College London,
London SW7 2AZ,
m.bridson@imperial.ac.uk
\medskip

KV: Mathematics Department, 555 Malott Hall, Cornell University,
Ithaca, NY 14850,
vogtmann@math.cornell.edu

 \end{document}